\documentclass[smallextended]{svjour3a}

\nonstopmode

\smartqed  
\usepackage{graphicx}
\usepackage{mathptmx}      

\journalname{Arnold Mathematical Journal}

\def\version{
September 28, 2019
}

\usepackage[usenames,dvipsnames]{color} 

\usepackage{latexsym,epsfig,bm}

\usepackage{bm}
\usepackage{upgreek}

\usepackage{mathrsfs}
\usepackage{times}

\definecolor{MyDarkBlue}{rgb}{0,0.08,0.45}
\definecolor{Pomegranate}{rgb}{0.6,0.1,0.15}
\definecolor{purple}{rgb}{0.6,0.1,0.15}
\usepackage[backref=none]{hyperref}
\hypersetup{pdfborder={0 0 0},
  colorlinks,
  urlcolor={MyDarkBlue},
  linkcolor={MyDarkBlue},
  citecolor={MyDarkBlue},
  breaklinks=true}
\usepackage{doi}

\usepackage{amssymb}
\usepackage{amsmath}

\providecommand{\eqref}[1]{{\rm (\ref{#1})}}

\newcommand{\unity}{\textrm{{\usefont{U}{fplmbb}{m}{n}1}}}

\DeclareSymbolFont{AMSb}{U}{msb}{m}{n}
\DeclareSymbolFontAlphabet{\mathbb}{AMSb}

\DeclareSymbolFont{EUR}{U}{eur}{m}{n}
\SetSymbolFont{EUR}{bold}{U}{eur}{b}{n}
\DeclareSymbolFontAlphabet{\eur}{EUR}

\DeclareSymbolFont{EUB}{U}{eur}{b}{n}
\SetSymbolFont{EUB}{bold}{U}{eur}{b}{n}
\DeclareSymbolFontAlphabet{\eub}{EUB}

\DeclareSymbolFont{EUS}{U}{eus}{m}{n}
\SetSymbolFont{EUS}{bold}{U}{eus}{b}{n}
\DeclareSymbolFontAlphabet{\eus}{EUS}

\newcommand{\jj}{\mathrm{i}}

\newcommand{\ch}{\,\mbox{\rm conv}\,}
\newcommand{\chomega}{\,\mbox{\rm conv}\sb{\!\!\sb\omega}\,}
\newcommand{\astomega}{\,\ast\sb{\!\!\sb\omega}\,}

\newcommand{\epi}{\mathop{\mathrm{epi}}}
\newcommand{\hyps}{\mathop{\mathrm{hyp}\sb{S}}}

\newcommand{\supp}{\mathop{\rm supp}}

\newcommand{\p}{\partial}
\newcommand{\at}[1]{\vert\sb{\sb{#1}}}

\def\R{\mathbb{R}}
\newcommand{\C}{\mathbb{C}}

\newcommand{\N}{\mathbb{N}}\newcommand{\Z}{\mathbb{Z}}

\newcommand{\abs}[1]{\vert #1 \vert}

\newcommand{\norm}[1]{\Vert #1 \Vert}
\newcommand{\sothat}{\,\,{\rm ;}\ \,}

\DeclareMathSymbol{\varGamma}{\mathord}{letters}{"00}
\DeclareMathSymbol{\varDelta}{\mathord}{letters}{"01}
\DeclareMathSymbol{\varTheta}{\mathord}{letters}{"02}
\DeclareMathSymbol{\varLambda}{\mathord}{letters}{"03}
\DeclareMathSymbol{\varXi}{\mathord}{letters}{"04}
\DeclareMathSymbol{\varPi}{\mathord}{letters}{"05}
\DeclareMathSymbol{\varSigma}{\mathord}{letters}{"06}
\DeclareMathSymbol{\varUpsilon}{\mathord}{letters}{"07}
\DeclareMathSymbol{\varPhi}{\mathord}{letters}{"08}
\DeclareMathSymbol{\varPsi}{\mathord}{letters}{"09}

\newtheorem{assumption}[lemma]{Assumption}

\renewenvironment{proof}[1][\proofname]{\medskip\noindent{\bfseries \it #1. }\ }{\qed\medskip}

\newcounter{step}

\makeatletter\@addtoreset{equation}{section}
\makeatletter\@addtoreset{lemma}{section}
\makeatother

\def\Mu{\mathrm{M}}
\def\Nu{\mathrm{N}}

\def\dist{\mathop{\rm dist}\nolimits}

\newcommand{\upsi}{u}

\renewcommand{\aa}[1]{\mathrm{A}\sb{{#1}}}
\newcommand{\bb}[1]{\mathrm{B}\sb{{#1}}}
\renewcommand{\bb}[1]{\mathrm{B}\sb{{#1}}}

\renewcommand{\aa}[1]{\eur{A}\sb{{#1}}}
\renewcommand{\bb}[1]{\eur{B}\sb{{#1}}}

\begin{document}

\title{
Solutions with compact time spectrum
to nonlinear Klein--Gordon and Schr\"odinger equations
and the Titchmarsh theorem for partial convolution
}

\titlerunning{Solutions with compact time spectrum
to NLKG and NLS
and the Titchmarsh theorem for partial convolution
}


\author{Andrew Comech}

\institute{
IITP, Moscow, Russia
\and
Texas A\&M University,
College Station, Texas, USA
\and
\email{comech@gmail.com}
}

\date{\version}

\maketitle

\begin{abstract}
We prove that
finite energy solutions
to the nonlinear Schr\"odinger equation
and nonlinear Klein--Gordon equation
which have the compact time spectrum
have to be one-frequency solitary waves.
The argument is based on the generalization of the Titchmarsh convolution theorem
to partial convolutions.
\end{abstract}


\bigskip
\bigskip

\hfill
{\it To Rafail Kalmanovich Gordin on the occasion of his 70th birthday
-- with love and admiration
}

\section{Introduction}

Let us consider nonlinear Schr\"odinger and
nonlinear Klein--Gordon equations,
\begin{eqnarray}\label{eqns}
\jj\p_t\upsi
=
-\Delta\upsi+\alpha(\abs{\upsi}^2)\upsi,
\qquad
-\p_t^2\upsi
=
-\Delta\upsi+m^2\upsi+\alpha(\abs{\upsi}^2)\upsi,
\qquad
\upsi(x,t)\in\C,
\qquad
x\in\R^n,
\quad
n\in\N,
\end{eqnarray}
where
$m>0$
and the nonlinearity is represented by a function
$\alpha\in C^1(\R)$, $\alpha(0)=0$.
These $\mathbf{U}(1)$-invariant equations
are well-known to admit solitary wave solutions of the form
\begin{eqnarray}\label{solitary-waves}
\upsi(x,t)=\phi(x)e^{-\jj\omega t},
\qquad
\omega\in\R,
\end{eqnarray}
with $\phi(x)$ decaying at infinity \cite{MR0454365,MR695535}.
Do these equations admit multifrequency solitary wave solutions
of the form $\sum_{j=1}^{N}\phi_j(x)e^{-\jj\omega_j t}$?
Indeed, such solutions have been found in similar systems;
see below for more details.
More generally, we would like to know
whether
besides one-frequency solitary waves
there are finite energy solutions with compact time-spectrum,
defined as follows.

\begin{definition}\label{def-spectrum}
Let $\upsi\in\mathscr{S}'(\R^n\times\R,\C)$,
and let
$\tilde\upsi(x,\omega)=\int\sb{\R} e^{\jj\omega t}\upsi(x,t)\,dt$
be its partial Fourier transform in time.
We say that the time spectrum of $\upsi$ is compact if
there is a finite interval
$I\subset\R$ such that
\[
\supp\tilde\upsi\subset\R^n\times I.
\]
\end{definition}

In the present article,
in Section~\ref{sect-kg},
we will prove that
in the nonlinear Schr\"odinger or Klein--Gordon equations
under certain assumptions on the nonlinearity (polynomial or some algebraic functions),
there are no finite energy solutions with compact time spectrum
except the one-frequency solitary waves of the form
\eqref{solitary-waves}.
See Theorem~\ref{theorem-one-frequency} below
for the precise formulation.
The approach is based on the form of the Titchmarsh convolution theorem
reformulated for partial convolutions;
see
Section~\ref{sect-titchmarsh}
and in particular Theorem~\ref{theorem-partial}.

\medskip

\noindent
{\bf Soliton resolution conjecture.\,}
This conjecture states that the long-time asymptotics
of any finite energy solution
to a nonlinear dispersive system
with $\mathbf{U}(1)$-symmetry
is given by a superposition of outgoing solitary waves
and an outgoing dispersive wave;
see
\cite{MR2032730,MR2275691,MR2304091,MR2308860,komech2016attractors}.
For the recent results
for the Schr\"{o}dinger and Klein--Gordon equations
with the critical power nonlinearity,
see \cite{duyckaerts2016} and the references therein.
Let us also mention the probabilistic approach \cite{MR3263670,MR3375596}.
One strategy to attack this problem was proposed in
\cite{MR2032730}:
one notices that any solution converges to \emph{radiationless solution},
the one that does not lose the energy any more.
Then one needs to complete the following two steps:

\begin{verse}
\it
1.
Prove that
any radiationless solution
has a compact time spectrum;

2.
Prove that any solution with compact time spectrum
has a time spectrum consisting of a single point,
and hence is a solitary wave:
\begin{eqnarray}\label{s-w-1}
\upsi(x,t)=
\phi(x)e^{-\jj\omega t},
\qquad
\omega\in\R,
\qquad
\phi\in H^1(\R^n,\C).
\end{eqnarray}
\end{verse}

Above, $H^1(\R^n)
=\{
u\in L^2(\R^n)
\sothat
\norm{u}_{H^1}^2
:=
\norm{u}_{L^2}^2+\norm{\nabla u}_{L^2}^2<\infty
\}$
is the standard Sobolev space of order one.

Thus, any finite energy solution converges
to a radiationless solution, which in turn is a solitary wave.
Both steps of the program were accomplished
for several models without translation invariance,
namely, for the Klein--Gordon equation
interacting with one oscillator
\cite{MR2032730,MR2308860},
\[
-\p_t^2\upsi=-\p_x^2\upsi+m^2\upsi
+\delta(x)\alpha(\abs{\upsi}^2)\upsi,
\qquad
\upsi(x,t)\in\C,
\quad
x\in\R,
\]
where $m>0$ and $\alpha(\tau)$ is a polynomial,
for several nonlinear oscillators
\cite{MR2579377},
for the Klein--Gordon and Dirac equations with the mean-field
self-interaction \cite{MR2526405,MR2745798}
(in any spatial dimension),
for the Klein--Gordon with the mean-field
self-interaction
at several points \cite{MR2902120},
and also for the Klein--Gordon equation in the discrete time-space
coupled to a nonlinear oscillator
\cite{MR3007724}.
In other words, in the models mentioned above,
\emph{the weak global attractor
is formed by solitary waves:}
any finite energy solution
converges to the solitary manifold,
\[
\mathbf{S}=\big\{
\phi\sb\omega(x)e^{-\jj\omega t}
\sothat
\ \omega\in\R
\big\}.
\]
The convergence is in the weak topology,
in weighted spaces
such as $H^1_{-s}(\R^n)=\{u\in H^1_{\mathrm{loc}}(\R^n)
\sothat \langle x\rangle^{-s} u\in H^1(\R^n)\}$, with $s>0$,
where $\langle x\rangle$ is [the operator of multiplication by]
the function $(1+x^2)^{1/2}$, $x\in\R^n$;
in this sense, we are talking about the \emph{weak attractor}.
The weight makes sure that we forget about the
excess energy, which is being carried away by the dispersive waves.
One then says that the convergence to the attractor is
caused by \emph{friction by dispersion}; this is the substitute
for the dissipation which is absent in a hamiltonian system.


In the present article we prove that,
under certain assumptions on the nonlinearity,
any solution with a compact time spectrum
is a single-frequency solitary wave.

\medskip

\noindent
{\bf Multifrequency solitary waves.\,}
If a particular model admits multifrequency solutions,
defined as exact localized solutions with several frequencies,
then they also belong to the attractor.
One can show that
multifrequency solitary waves exist in the Klein--Gordon equation
with the mean-field self-interaction \cite{MR2526405}
and with several nonlinear oscillators  \cite{MR2579377}.
Bi-frequency solitary waves exist in systems of nonlinear
Schr\"odinger equations \cite{PhysRevA.86.053809}
(the vector case may admit solutions with several harmonics
when the nonlinearity does not produce higher harmonics
due to cancellations, which are absent in the scalar case).
In a similar fashion, bi-frequency solitary waves
exist in the Soler model and Dirac--Klein--Gordon model with Yukawa self-interaction
\cite{MR3842864}:
\[
\psi(x,t)=\phi(x)e^{-\jj\omega t}+\chi(x)e^{\jj\omega t},
\qquad
\mbox{for particular}
\ \ \phi,\,\chi\in H^1(\R^n,\C^N).
\]
Sometimes one may place some restriction on the parameters of the
problem (such as the spacings between the nonlinear
oscillators in \cite{MR2579377}) to ensure that
multifrequency solutions would be absent.

In \cite{MR3007724},
based on the Titchmarsh theorem for distributions on the circle \cite{MR3087829},
it was shown that
the global attractor of the Klein--Gordon equation in discrete time-space
coupled with a nonlinear oscillator,
besides usual one-frequency solitary waves $\phi e^{-\jj\omega T}$,
could also contain two- and four-frequency solutions:
\[
\phi e^{-\jj\omega T}
+
\chi e^{-\jj(\omega+\pi)T},
\qquad
\qquad
\phi_1 e^{-\jj\omega_1 T}
+\phi_2 e^{-\jj\omega_2 T}
+\chi_1 e^{-\jj(\omega_1+\pi)T}
+\chi_2 e^{-\jj(\omega_2+\pi)T},
\]
where $T\in\Z$ is the discrete time
and $\phi,\,\chi,\,\ldots\in l^2(\Z^n)$ are particular functions
of the discrete spatial variable $X\in\Z^n$,
and indeed examples of such solutions were given.

According to Theorem~\ref{theorem-one-frequency} (see below),
the nonlinear Schr\"odinger and Klein--Gordon equations
with a certain class of nonlinearities
do not admit multifrequency solitary wave solutions.

\medskip

\noindent
{\bf Breathers.\,}
Let us contrast our results to the existence of breathers,
which are exact periodic solutions
in the context of completely integrable systems.
For example, the completely integrable sine--Gordon equation
\begin{eqnarray}\label{sine-gordon}
-\p_t^2 u=-\p_x^2 u+\sin u,
\qquad
u(x,t)\in\R,
\quad
x\in\R,
\end{eqnarray}
admits
solutions of the following form
\cite{PhysRevLett.30.1262}:
\[
u(x,t)=4\arctan
\left(
\frac{\sqrt{1-\omega^2}\cos(\omega t)}
{\omega\cosh(\sqrt{1-\omega^2}\,x)}
\right),
\qquad
\omega\in(-1,1),
\]
which are exponentially localized in space and are periodic in time.
Note that the time spectrum of this solution is unbounded,
and moreover the nonlinearity
in \eqref{sine-gordon}
is not of algebraic type;
thus, this solution does not contradict
our statement on the absence of nontrivial compact spectrum solutions
(other than one-frequency solitary waves)
to the nonlinear Klein--Gordon equation with certain algebraic nonlinearities.

Similarly, the cubic nonlinear Schr\"odinger equation
\[
\jj\p_t u=-\p_x^2 u-2\abs{u}^2 u,
\qquad
u(x,t)\in\C,
\quad
x\in\R,
\]
admits exact solutions \cite{MR915545}
such as the following one:
\[
u(x,t)=
\frac
{\cos x+\jj\sqrt{2}\sinh t}
{\sqrt{2}\cosh t-\cos x}
e^{\jj t}.
\]
We notice that the frequency spectrum of this solution
is not compact; moreover, this solution has an infinite $L^2$-norm
and energy.
For more examples of such solutions, see \cite{MR915545}.

\medskip

\noindent
{\bf Convergence of small initial data to one-frequency solitary waves.\,}
Let us mention the results on convergence
of small solutions to (one-frequency) solitary waves,
particularly in the context of the nonlinear Schr\"{o}dinger equation:
in other words,
the attractor of small solutions is formed
by small amplitude solitary waves.
See in particular
\cite{tsai2002relaxation,soffer2004,cuccagna2015small,cuccagna2016.32,cuccagna2016.1332}.

\section{Titchmarsh theorem for partial convolution}
\label{sect-titchmarsh}

The original formulation of the Titchmarsh convolution theorem
\cite{titchmarsh1926zeros} is as follows:
\begin{verse}
{\it 
\noindent
If $\phi(t)$ and $\psi(t)$ are integrable functions, such that
\mbox{$
\displaystyle\int\sb{0}\sp{x}\phi(t)\psi(x-t)\,dt=0$}
almost everywhere in the interval $0<x<\kappa$, then $\phi(t)=0$ almost
everywhere in $(0,\lambda)$,
and $\psi(t)=0$ almost everywhere in $(0,\mu)$, where
$\lambda+\mu\ge\kappa$.
}
\end{verse}

Above, $\lambda$ and $\mu$ are some particular values $\ge 0$.
An equivalent reformulation is that
$\inf\supp\phi\ast\psi
=\inf\supp \phi+\inf\supp\psi$,
for any
$\phi,\,\psi\in\mathscr{E}'(\R)$,
where $\mathscr{E}'(\R)$
is the space of distributions with compact support
(dual to the space $\mathscr{E}(\R)$
which is $C\sp\infty(\R)$
with the topology defined by the seminorms
$\sup_{\omega\in K} |f^{(k)}(\omega)|$,
with $k\in\mathbb{N}_0$ and $K$ a compact subset of $\mathbb{R}$).
A higher dimensional generalization
can be stated in terms of the convex hulls of the supports
of distributions \cite{MR0043254}:

\begin{theorem}[Titchmarsh Convolution Theorem
\cite{MR0043254}]
For $f,g\in\mathscr{E}'(\R^n)$,
\begin{equation}
\ch\supp f\ast g=\ch\supp f+\ch\supp g.
\end{equation}
\end{theorem}

Above,
$\ch$ denotes the convex hull of a set.

We need a version of this theorem for a partial convolution
with respect to only a subset of variables.

\subsection{Maximal lower semicontinuous function
and minimal upper semicontinuous function}

\begin{lemma}\label{lemma-upper}
Let $n\ge 1$.
For any function $\mu:\R^n\to\R$
there is a maximal lower semicontinuous function
on $\R^n$
which does not exceed $\mu$;
we will denote this function
by $\mu^L(x)$.
Similarly,
there is a minimal upper semicontinuous function
on $\R^n$
which is not exceeded by $\mu$;
we will denote this function
by $\mu^U(x)$.
For any $\mu,\,\nu:\R^n\to\R$,
one has
\begin{eqnarray}\label{l-u}
&
\mu^L\le \mu\le\mu^U,
\\
\label{not-in-general}
&
(\mu+\nu)^L\ge\mu^L+\nu^L,
\qquad
(\mu+\nu)^U\le\mu^U+\nu^U.
\end{eqnarray}
\end{lemma}

\begin{proof}
The function is lower semicontinuous
if and only if its \emph{epigraph}
(the set of points lying on or above its graph),
$
\epi\mu
=\{(x,y)\in \R^n\times\R\sothat y\ge \mu(x)\},
$
is closed,
or, equivalently, if and only if its \emph{strict epigraph},
\[
\hyps(\mu)
=\{(x,y)\in \R^n\times\R\sothat y<\mu(x)\},
\]
is open.
For a function $\mu:\R^n\to\R$ let us consider
the complement to its epigraph, the strict \emph{hypograph},
\[
\hyps(\mu)
=\{(x,y)\in \R^n\times\R\sothat y<\mu(x)\}.
\]
Let $\mu\sb\alpha:\R^n\to\R$,
$\alpha\in I$,
be a subset of the set of lower semicontinuous functions.
Then
\[
\hyps\Big(\sup\sb{\alpha\in I}\mu\sb\alpha\Big)
=
\mathop{\cup}\limits\sb{\alpha\in I}
\hyps(\mu\sb\alpha)
\]
is open
(as a union of any collection of open sets),
hence
$\mu^L:=\sup\sb{\alpha\in I}\mu\sb\alpha$
is lower semicontinuous.

The inequalities
\eqref{l-u} and \eqref{not-in-general}
readily follow from the definition of $\mu^L$ and $\mu^U$.
\end{proof}

\begin{remark}\label{remark-example}
The example of upper semicontinuous functions
\[
\mu=\begin{cases}1,&x\le 0
\\0,&x> 0\end{cases},
\qquad
\nu=\begin{cases}0,&x<0
\\1,&x\ge 0\end{cases},
\qquad
\mu+\nu=\begin{cases}1,&x\ne 0
\\2,&x=0\end{cases},
\]
with
$
\mu^L=\begin{cases}1,&x< 0
\\0,&x\ge 0\end{cases},
$
$
\quad \nu^L=\begin{cases}0,&x\le 0
\\1,&x>0\end{cases},
$
$
\quad (\mu+\nu)^L\equiv 1,
$
shows that the strict inequalities
in \eqref{not-in-general} are possible
(for $\R$-valued functions
we say that $f<g$ if there is at least one point $x$
in their domains such that $f(x)<g(x)$).
\end{remark}




We recall that the space of distributions
$\mathscr{D}'(\R^n)$ is defined as the dual to
$\mathscr{D}(\R^n)=C^\infty_{\mathrm{comp}}(\R^n)$,
while
$\mathscr{E}'(\R^n)$ is the space of distributions
with compact support
(the dual to $C^\infty(\R^n)$).

\begin{definition}
Let
$f\in\mathscr{D}'(\R^n\times\R)$.
We define the functions $\aa{f}$ and $\bb{f}$ by
\begin{eqnarray}
\nonumber
&&
\aa{f}:\;\R^n\to\R\sqcup\{\pm\infty\},
\qquad
x\mapsto
\inf\big\{
\omega\in\R\sothat
(x,\omega)\in\supp f
\big\};
\\[1ex]
\nonumber
&&
\bb{f}:\;\R^n\to\R\sqcup\{\pm\infty\},
\qquad
x\mapsto
\sup\big\{
\omega\in\R\sothat
(x,\omega)\in\supp f
\big\}.
\end{eqnarray}
\end{definition}

It follows that
$\aa{f}$ is lower semicontinuous,
while
$\bb{f}$ is upper semicontinuous:
\[
\aa{f}=\aa{f}^L,
\qquad
\bb{f}=\bb{f}^U.
\]

\begin{definition}\label{def-sigma}
Let $f\in\mathscr{D}'(\R^n\times\R)$.
We define
$\Sigma[f]$ to be the projection of
$\supp f\subset\R^n\times\R$ onto the first factor:
\[
\Sigma[f]
=\big\{x\in\R^n\sothat
(\{x\}\times\R)\cap\supp f\ne\emptyset
\big\}\subset\R^n.
\]
\end{definition}

Thus, one has
\[
x\not\in\Sigma[f]
\quad\Leftrightarrow\quad
\aa{f}(x)=+\infty
\quad\Leftrightarrow\quad
\bb{f}(x)=-\infty.
\]

\begin{lemma}\label{lemma-sigma-closed}
Let $f\in \mathscr{D}'(\R^n\times\R)$.
If there is a finite interval $I\subset\R$
such that $\supp f\subset\R^n\times I$,
then the set $\Sigma[f]\subset\R^n$ is closed.
\end{lemma}

\begin{remark}
$\Sigma[f]$ is not necessarily closed for
$f\in \mathscr{D}'(\R^n\times\R)$.
\end{remark}

\begin{lemma}\label{lemma-a-a-b}
For any distribution $f\in\mathscr{D}'(\R^n\times\R)$,
one has
\begin{eqnarray}\label{a-a-b}
\aa{f}(x)\le \aa{f}^U(x)\le \bb{f}(x),
\qquad
\aa{f}(x)\le \bb{f}^L(x)\le \bb{f}(x),
\qquad
\forall x\in\Sigma[f];
\end{eqnarray}
\begin{eqnarray}\label{a-u-l-a}
(\aa{f}^U)^L\ge \aa{f},
\qquad
(\bb{f}^L)^U\le \bb{f}.
\end{eqnarray}
\end{lemma}

\begin{proof}
Note that $\bb{f}$ is upper semicontinuous with
$\bb{f}(x)\ge \aa{f}(x)$
for all $x\in\Sigma[f]$,
while $\aa{f}^U$ is the
\emph{smallest} upper semicontinuous function
which is not smaller than $\aa{f}$
(cf. Lemma~\ref{lemma-upper});
thus,
$\aa{f}\le \aa{f}^U\le \bb{f}$.
The second relation in \eqref{a-a-b} is proved similarly.

For the relations \eqref{a-u-l-a},
one can see that for any lower semicontinuous function
$a:\;\R^n\to\R$ one has
$(a^U)^L\ge a$
(since $a$ is a lower semicontinuous function
which is not larger than $a^U$),
and similarly for any upper semicontinuous function $b:\R^n\to \R$,
one has $(b^L)^U\le b$.
\end{proof}

\begin{remark}
In \eqref{a-a-b},
$\aa{f}^U$ is not necessarily smaller than
$\bb{f}^L$;
it suffices to consider the example
$f(x,\omega)=\theta(-x)\delta_{-1}(\omega)
+
\theta(x)\delta_{1}(\omega)$,
with $x,\,\omega\in\R$.
Also, the inequalities in \eqref{a-u-l-a}
are not necessarily strict,
as the example
$f(x,\omega)=\delta(\omega)+\delta(x)\unity_{[-1,1]}(\omega)$
shows
(in detail, $\bb{f}(0)=1$, $\bb{f}(x)=0$ for $x\ne 0$;
$\aa{f}(0)=-1$, $\bb{f}(x)=0$ for $x\ne 0$;
$\bb{f}^L\equiv 0\equiv \aa{f}^U$,
$(\bb{f}^L)^U\equiv 0\equiv (\aa{f}^U)^L$).
\end{remark}


For $f,\,g\in C^\infty\sb{\mathrm{comp}}(\R^n\times\R)$,
we define the partial convolution
\begin{eqnarray}\label{partial-convolution}
&&
\nonumber
\astomega:\;
C^\infty\sb{\mathrm{comp}}(\R^n\times\R)
\times
C^\infty\sb{\mathrm{comp}}(\R^n\times\R)
\to
C^\infty\sb{\mathrm{comp}}(\R^n\times\R),
\\[1ex]
&&
(f\astomega g)(x,\omega)
=\int\sb{\R}f(x,\omega-\tau)g(x,\tau)\,d\tau,
\qquad
(x,\omega)\in\R^n\times\R.
\end{eqnarray}
This operation can be continuously extended to
$f,\,g\in \mathscr{E}'(\R,L^2(\R^n))$:
\[
\astomega:\;
\mathscr{E}'(\R,L^2(\R^n))
\times
\mathscr{E}'(\R,L^2(\R^n))
\to
\mathscr{E}'(\R,L^1(\R^n)).
\]
Indeed, let $f,\,g,\,\varphi\in C^\infty\sb{\mathrm{comp}}(\R^n\times\R)$.
Then
\[
\big\langle
f\astomega g,\varphi
\big\rangle
=
\Big\langle
\int_\R f(x,\omega-\tau)g(x,\tau)\,d\tau,\,\varphi(x,\omega)
\Big\rangle
=
\int_{\R}
\Big(
g(x,\tau)\,
\int_{\R^n\times\R} f(x,\omega-\tau)\varphi(x,\omega)\,dx\,d\omega
\Big)\,d\tau,
\]
where $\langle\,\cdot\,,\,\cdot\,\rangle$
refers to the pairing
of $L^2(\R^n\times\R)$-functions.
The integral
$
\int_{\R^n\times\R} f(x,\omega-\tau)\varphi(x,\omega)\,dx\,d\omega
$
makes sense
for $f\in\mathscr{E}'(\R,L^2(\R^n))$
and
$\varphi\in\mathscr{E}(\R,L^\infty(\R^n))$,
defining an element from $\mathscr{E}(\R,L^2(\R^n))$,
which could then be coupled with
$g\in\mathscr{E}'(\R,L^2(\R^n))$.
Let us mention that
for (complex) Banach spaces $A,\,B$
and the space of bounded linear maps $\mathscr{B}(A,B)$,
the space of $\mathscr{B}(A,B)$-valued distributions
$\mathscr{D}'(\R,\mathscr{B}(A,B))$
is defined as the space of bounded linear maps
from $\mathscr{D}(\R,A)$
($A$-valued test functions)
to $B$,
and similarly for $\mathscr{B}(A,B)$-valued tempered distributions $\mathscr{E}'$;
for the general theory of Banach-space-valued distributions,
see \cite[Chapter 3]{zemanian1972}.

\medskip

For any
$f,\,g\in\mathscr{E}'(\R,L^2(\R^n))$
there are immediate relations
\begin{eqnarray}\label{immediate}
\aa{f\astomega g}\ge \aa{f}+\aa{g},
\qquad
\bb{f\astomega g}\le \bb{f}+\bb{g}.
\end{eqnarray}

We will show that the relations \eqref{immediate} are
equalities, in the appropriate sense.


\begin{theorem}[Titchmarsh theorem for partial convolution]~
\label{theorem-partial}
Let $f,\,g\in \mathscr{E}'(\R,L^2\sb{\mathrm{loc}}(\R^n))$.
Then
\[
\aa{f\astomega g}
=
\big(\aa{f}^U+\aa{g}\big)^L
=\big(\aa{f}+\aa{g}^U\big)^L,
\qquad
\bb{f\astomega g}
=
\big(\bb{f}^L+\bb{g}\big)^U
=\big(\bb{f}+\bb{g}^L\big)^U.
\]
\end{theorem}

\begin{remark}
Let us prove a similar statement
for elements from the space $C(\R^n,\mathscr{E}'(\R))$,
defined as the space of functions
$F:\,\R^n\to\mathscr{E}'(\R)$
which satisfy
$\lim\sb{x\to x_0}F(x)=F(x_0)$
for any $x_0\in\R^n$,
with the convergence in the topology of $\mathscr{E}'(\R)$.
For $f,\,g\in C(\R^n,\mathscr{E}'(\R))$,
since $f$ and $g$ depend continuously on $x$,
the Titchmarsh convolution theorem
can be applied pointwise in $x$,
yielding
\begin{eqnarray}\label{standard-Titchmarsh}
\inf\supp (f\astomega g)(x,\cdot)
=
\inf\supp f(x,\cdot)
+
\inf\supp g(x,\cdot),
\qquad
\forall x\in\R^n,
\end{eqnarray}
and similarly for $\sup$.
Let $f\in C(\R^n,\mathscr{E}'(\R))$
and let $\rho\in\mathscr{D}(\R)$.
If $\mathcal{O}\subset\R^n$ is an open set such that
$\langle \rho,f(x,\cdot)\rangle=0$
for all $x\in \mathcal{O}$,
then, by continuity of $f$ in $x$,
one also has $\langle \rho,f(x,\cdot)\rangle=0$
for all $x$ from the closure of $\mathcal{O}$.
Therefore, given an open set $\Omega\subset\R$,
if $\Omega\cap\supp f(x,\cdot)=\emptyset$ for $x\in \mathcal{O}\subset\R^n$,
then
$\Omega\cap\supp f(x,\cdot)=\emptyset$ for $x$ from the closure of $\mathcal{O}$;
it then follows that
\[
\aa{f}^U(x)=\inf\supp f(x,\cdot),
\qquad
\bb{f}^L(x)=\sup\supp f(x,\cdot)
\qquad
\mbox{for any}\ x
\ \mbox{such that}\ \;f(x,\cdot)\not\equiv 0.
\]
Applying the above to each of the terms
in \eqref{standard-Titchmarsh}
(and similarly for $\sup\supp$)
leads to the relations
\[
\aa{f\astomega g}^U(x)
=
\aa{f}^U(x)
+
\aa{g}^U(x),
\qquad
\bb{f\astomega g}^L(x)
=
\bb{f}^L(x)
+
\bb{g}^L(x)
\]
which are similar to the relations stated in Theorem~\ref{theorem-partial}.
\end{remark}







\subsection{Convex hulls and partial convolution theorem in higher dimensions}

Let us give a higher dimensional version
of the partial convolution theorem
in terms of convex hulls, following \cite{MR0043254}.
Let $n,\,m\ge 1$.
For any set-valued map $\Mu:\R^n\to\{\mbox{closed subsets of $\R^m$}\}$
there is a maximal inner semicontinuous set-valued map
$\R^n\to\{\mbox{closed subsets of $\R^m$}\}$
which does not exceed $\Mu$;
we denote this map by
\[
\Mu^L(x)
=
\lim\sb{\epsilon\to 0}
\bigcap\sb{y\in\mathbb{B}_\epsilon(x)}\Mu(y),
\qquad\forall x\in\R^n.
\]
Note that for each $x\in\R^n$, the set $\Mu^L(x)\subset\R^m$
is closed (as an intersection of an arbitrary number of closed sets).
Similarly,
there is a minimal outer semicontinuous set-valued map
$\R^n\to\{\mbox{closed subsets of $\R^m$}\}$
which is not exceeded by $\Mu$;
we denote this map by
\[
\Mu^U(x)
=
\lim\sb{\epsilon\to 0}
\bigcup\sb{y\in\mathbb{B}_\epsilon(x)}\Mu(y),
\qquad\forall x\in\R^n.
\]
(Note that for each $x\in\R^n$,
the set
$\Mu^U(x)\subset\R^m$ is closed:
if $\omega_j\in \Mu^U(x)$
converges to some $\omega_*\in\R^m$ as $j\to\infty$, then there are
sequences $\omega_{j,N}\in \Mu(y_N)$ with $|x-y_N|<1/N$, $N\in\N$
such that, for each $j\in\N$,
$\omega_{j,N}\to \omega_j$ as $N\to\infty$,
but  then one can choose a diagonal subsequence $\omega_{j_r,N_r}$
converging to $\omega_*$.
Thus,
$\omega_*\in \Mu^U(x)$, so $\Mu^U(x)$ is closed.)
Thus,
\[
\Mu^L(x)\subset \Mu(x)\subset\Mu^U(x),
\qquad
\forall x\in\R^n.
\]

The following lemma is an immediate generalization of
Lemma~\ref{lemma-upper}.

\begin{lemma}\label{lemma-upper-set}
For any $\Mu,\,\Nu:\R^n\to\{\mbox{closed subsets of $\R^m$}\}$,
one has
\begin{eqnarray}\label{l-u-set}
\label{not-in-general-set}
(\Mu+\Nu)^L(x)\supset\Mu^L(x)+\Nu^L(x),
\qquad
(\Mu+\Nu)^U(x)\subset\Mu^U(x)+\Nu^U(x),
\qquad
\forall x\in\R^n.
\end{eqnarray}
\end{lemma}

Above, 
the sum of two subsets $A,\,B\subset\R^m$
is defined by
$A+B=\{a+b\in\R^m\sothat a\in A,\,b\in B\}\subset\R^m$.

We recall that, given a set $S\subset\R^n$,
then $\ch S$
denotes its convex hull.
For a set $S\subset\R^n\times\R^m$,
with $m,\,n\in\N$,
let us define
$\chomega S$
as a map from $\R^n$ to convex subsets of $\R^m$ by
\[
\chomega S:
\;x\mapsto
\ch\big(S\cap\big(\{x\}\times\R^m\big)\big)
\subset\R^m.
\]
If $S$ is closed, this map
is outer semicontinuous.

For a closed subset $S\subset\R^n\times\R^m$,
we define
\[
(\chomega S)^L:\;
\R^n
\to
\big\{\mbox{closed subsets of $\R^n\times\R^m$}\big\}
\]
as the largest inner semicontinuous map
from $\R^n$ to closed convex subsets of $\R^m$
which satisfies
\[
(\chomega S)^L(x)\subset(\chomega S)(x),
\qquad\forall x\in\R^n.
\]

\begin{remark}
For $f\in\mathscr{D}'(\R^n\times\R)$,
there is an obvious relation
\[
(\chomega\supp f)^L(x)=\big[\aa{f}^U(x),\bb{f}^L(x)\big]
\subset
(\chomega\supp f)(x)=\big[\aa{f}(x),\bb{f}(x)\big],
\qquad
\forall x\in\Sigma[f].
\]
\end{remark}

\begin{theorem}[Titchmarsh theorem for partial convolution: convex hulls]
\label{theorem-partial-2}
Let $f,\,g\in
\mathscr{E}'(\R^m,L^2\sb{\mathrm{loc}}(\R^n))$.
Then
\[
\chomega\supp f\astomega g
=
\Big(
(\chomega\supp f)^L
+
\chomega\supp g
\Big)^U
=
\Big(
\chomega\supp f
+
(\chomega\supp g)^L
\Big)^U.
\]
\end{theorem}

The proof of Theorem~\ref{theorem-partial-2}
follows the same lines
as that of Theorem~\ref{theorem-partial}
(using the language of \cite{MR0043254}).

\subsection{Proof of partial convolution theorem for $f\astomega f$
for $f\in \mathscr{D}(\R,L^2_{\mathrm{loc}}(\R^n))$}

Following \cite[Proof of Theorem 4.3.3]{MR717035},
we first prove the theorem for
$f,\,g\in \mathscr{D}(\R,L^2\sb{\mathrm{loc}}(\R^n))$.
To consider the case $f=g$,
we need the two lemmata,
which are the immediate adaptations of \cite[Lemmata 4.3.4, 4.3.5]{MR717035}.

\begin{lemma}\label{lemma-4.3.4}
For
$f\in \mathscr{D}(\R,L\sb{\mathrm{comp}}^2(\R^n))$,
one has
\[
\norm{f\astomega f\sp\sharp}_{L^2(\R^n\times\R)}
=
\norm{f\astomega f}_{L^2(\R^n\times\R)},
\]
where
\begin{eqnarray}\label{def-sharp}
f\sp\sharp(x,\omega)=\overline{f(x,-\omega)}.
\end{eqnarray}
\end{lemma}

\begin{lemma}\label{lemma-4.3.5}
For any finite open interval $\Omega\subset\R$,
there is $C<\infty$ such that
\[
\norm{f}_{L^\infty(\R,L^2(\R^n))}
\le
C\norm{\p_\omega^2 f}_{L^2(\R^n\times \Omega)},
\qquad
\forall f\in
\mathscr{D}(\Omega,L^2_{\mathrm{comp}}(\R^n)).
\]
\end{lemma}

Now we can give the proof for the case $f=g
\in \mathscr{D}(\R,L^2\sb{\mathrm{loc}}(\R^n))$.

\begin{lemma}\label{lemma-f-f-smooth}
Let $f\in \mathscr{D}(\R,L^2\sb{\mathrm{loc}}(\R^n))$.
There are the relations
\ $
\aa{f\astomega f}=2\aa{f},
$
\ $
\bb{f\astomega f}=2\bb{f}.
$
\end{lemma}


\begin{proof}
For any open set $\mathcal{O}\subset\R^n$,
one has:
\[
\norm{f}^2_{L^2(\mathcal{O}\times \Omega)}
=
\norm{f\astomega f\sp\sharp(\cdot,0)}_{L^2(\mathcal{O})}
\le
\norm{f\astomega f\sp\sharp}_{L^2(\mathcal{O},L^\infty(\Omega))}
\]
\[
\le
C\norm{\p_\omega^2(f\astomega f\sp\sharp)}_{L^2(\mathcal{O}\times \Omega)}
=
C\norm{\p_\omega f\astomega \p_\omega f\sp\sharp)}_{L^2(\mathcal{O}\times \Omega)}
=
C\norm{\p_\omega f\astomega \p_\omega f)}_{L^2(\mathcal{O}\times \Omega)};
\]
in the second line,
we applied Lemma~\ref{lemma-4.3.5}
and then Lemma~\ref{lemma-4.3.4}.
Applying the above inequality
to $f_\xi(x,\omega)=f(x,\omega)e^{\omega\xi}$,
we arrive at the inequality
\begin{eqnarray}\label{u-u-u}
\norm{f_\xi}_{L^2(\mathcal{O}\times \Omega)}^2
\le
C\norm{\p_\omega^2 (f_\xi\astomega f_\xi)}_{L^2(\mathcal{O}\times \Omega)}.
\end{eqnarray}
This inequality is satisfied for arbitrarily large $\abs{\xi}$,
while
$
f_\xi\astomega f_\xi(x,\omega)
=
e^{\omega\xi}(f\astomega f)(x,\omega)
$
for a given function $f$;
hence
twice the support  of the integrand in the left-hand side
of \eqref{u-u-u}
is contained in
$\mathop{\cup}\limits_{x\in \mathcal{O}}
\ch\big((\{x\}\times \Omega)\cap\supp f\astomega f\big)$.
Sending $\mathcal{O}\to\{x\}$,
we conclude that
$
2 \bb{f}(x)\le \bb{f\astomega f}(x),
$
for all $x\in\R^n$.
We conclude that
\begin{eqnarray}\label{2-b-b-u-u}
2 \bb{f}(x)\le \bb{f\astomega f}(x),
\qquad
\forall x\in\R^n.
\end{eqnarray}
Due to an immediate inequality
$2 \bb{f}(x)\ge \bb{f\astomega f}(x)$
which follows from the definition
\eqref{partial-convolution},
one has
$2 \bb{f}(x)=\bb{f\astomega f}(x)$.
Similarly,
$2 \aa{f}(x)=\aa{f\astomega f}(x)$.
\end{proof}

\subsection{Proof of partial convolution theorem
for $f,\,g\in \mathscr{D}(\R,L^2_{\mathrm{loc}}(\R^n))$}
\label{sect-f-g-d}

\begin{lemma}\label{lemma-bd-c}
Let $f,\,g\in 
\mathscr{D}(\R,L^2\sb{\mathrm{loc}}(\R^n))$.
Then, for any polynomials $\alpha(\omega)$ and $\beta(\omega)$,
\[
\aa{(\alpha f)\astomega(\beta g)}\ge \aa{f\astomega g},
\qquad
\bb{(\alpha f)\astomega(\beta g)}\le \bb{f\astomega g}.
\]
\end{lemma}

\begin{proof}
We closely follow the argument from \cite[Proof of Theorem 4.3.3]{MR717035}.
It suffices to prove the second inequality,
and only for the polynomials
$\alpha(\omega)=\omega$, $\beta(\omega)=1$.
Denote
\begin{equation}\label{def-f-g-n}
f\sb n(x,\omega)=\omega^n f(x,\omega),
\qquad
g\sb n(x,\omega)=\omega^n g(x,\omega),
\qquad
B\sb{m n}(x):=\bb{f\sb m\astomega g\sb n}(x).
\end{equation}
Let us assume that, contrary to the statement of the Lemma,
there is $x\in\R^n$ such that
\begin{equation}\label{ass-a00}
\bb{f\sb 1\astomega g}(x)>\bb{f\astomega g}(x);
\end{equation}
from now on, all the quantities are evaluated at this particular value of $x$.
The inequality \eqref{ass-a00} can be rewritten as
\begin{equation}\label{ass-a00-1}
B\sb{10}-B\sb{00}>0.
\end{equation}
Due to the relation
$\displaystyle
\omega (f\astomega g)(\omega)=(f\sb 1\astomega g)(\omega)+(f\astomega g\sb 1)(\omega),
$
we have:
\begin{equation}\label{a00-fg-fg}
\bb{f\sb 1\astomega g+f\astomega g\sb 1}
=\bb{\omega(f\astomega g)(\omega)}
\le \bb{f\astomega g}=B\sb{00}.
\end{equation}
It follows that
\[
\bb{f\sb 1\astomega g\astomega f\sb 1\astomega g+f\sb 1\astomega g\astomega f\astomega g\sb 1}
\leq
\bb{f\sb 1\astomega g} +\bb{f\sb 1\astomega g+f\astomega g\sb 1}
\leq B\sb{10}+B\sb{00}.
\]
If we had
$\bb{f\sb 1\astomega g\astomega f\sb 1\astomega g}
\ne \bb{f\sb 1\astomega g\astomega f\astomega g\sb 1}$,
then both these quantities would be smaller than or equal to
$B\sb{10}+B\sb{00}$.
By Lemma~\ref{lemma-f-f-smooth}
and \eqref{ass-a00-1},
this would lead to
$\bb{f\sb 1\astomega g}\le(B\sb{10}+B\sb{00})/2<B\sb{10}$,
contradicting \eqref{def-f-g-n}.
Thus,
$
\bb{f\sb 1\astomega g\astomega f\sb 1\astomega g}
=
\bb{f\sb 1\astomega g\astomega f\astomega g\sb 1},
$
leading to
\begin{equation}\label{fg-fgfg-fg}
\bb{f\sb 1\astomega g\astomega f\sb 1\astomega g}
=\bb{f\sb 1\astomega g\astomega f\astomega g\sb 1}
\le
\bb{f\astomega g}+\bb{f\sb 1\astomega g\sb 1}.
\end{equation}
By Lemma~\ref{lemma-f-f-smooth},
$\bb{f\sb 1\astomega g\astomega f\sb 1\astomega }
=2 \bb{f\sb 1\astomega g}$;
then \eqref{fg-fgfg-fg} could be rewritten as
\begin{equation}\label{a0-a0}
2\bb{f\sb 1\astomega g}
\le
\bb{f\astomega g}+\bb{f\sb 1\astomega g\sb 1}.
\end{equation}
This gives
\begin{equation}\label{a10-a11}
B\sb{11}-B\sb{10}\ge B\sb{10}-B\sb{00}>0.
\end{equation}
In the last inequality, we took into account \eqref{ass-a00-1}.
The inequalities \eqref{a10-a11} imply that
\begin{equation}\label{ass-a11}
\bb{f\sb 1\astomega g\sb 1}>\bb{f\sb 1\astomega g}.
\end{equation}
Just as we derived \eqref{a0-a0} from \eqref{ass-a00},
we could use \eqref{ass-a11}
to derive
\begin{equation}\label{a1-a1}
2\bb{f\sb 1\astomega g\sb 1}
\le
\bb{f\sb 1\astomega g}+\bb{f\sb 2\astomega g\sb 1}.
\end{equation}
The inequality
\eqref{a1-a1}
could be written as
$B\sb{21}-B\sb{11}\geq B\sb{11}-B\sb{10}$,
and, together with \eqref{a10-a11}, this yields
\[
B\sb{21}-B\sb{11}\geq B\sb{11}-B\sb{10}
\geq B\sb{10}-B\sb{00}>0.
\]
Proceeding by induction,
we prove that
$
B\sb{32}-B\sb{22}\geq B\sb{22}-B\sb{21}
\geq
B\sb{21}-B\sb{11}\geq B\sb{11}-B\sb{10}
\geq B\sb{10}-B\sb{00}>0,
$
hence
\begin{equation}\label{a-n}
B\sb{n n}\ge B\sb{00}+2n(B\sb{10}-B\sb{00}).
\end{equation}
At the same time,
since $\bb{f\sb n}\le \bb{f}$,
$\bb{g\sb n}\le \bb{g}$,
we know that
$
\bb{f\sb n\astomega g\sb n}
\le
\bb{f\sb n}+\bb{g\sb n}
\le \bb{f}+\bb{g}.
$
This would be in contradiction to \eqref{a-n}.
Hence, \eqref{ass-a00} is not true.
This finishes the proof of the lemma.
\end{proof}

\begin{proof}[Proof of Theorem~\ref{theorem-partial}
for $f,\,g\in \mathscr{D}(\R,L^2_{\mathrm{loc}}(\R^n))$]
Now we complete the proof of the Titchmarsh
theorem for $f\astomega g$.
For our convenience, we assume that
$\supp f\subset\R^n\times[1,+\infty)$
and
$\supp g\subset\R^n\times[1,+\infty)$.

Fix $x\in\R^n$.
Let $\epsilon\in(0,1)$.
Due to lower semicontinuity of $\aa{f\astomega g}$,
for any $\omega_0\in\big(\aa{f\astomega g}(x)-\epsilon,\aa{f\astomega g}(x)\big)$,
there is a nonempty open neighborhood
$\mathcal{O}\subset\mathbb{B}_\epsilon(x)$,
$\mathcal{O}\ni x$,
such that
$\omega_0<\aa{f\astomega g}(y)$
for all $y\in \mathcal{O}$.
This implies that
\begin{equation}\label{int-f-g-n}
\int_{\mathcal{O}} \rho(y)
\int\sb{0}\sp\omega f(y,\omega-\tau)g(y,\tau)\,d\tau\,dy=0,
\qquad\forall\omega\in(0,\omega_0),
\qquad\forall \rho\in C^\infty\sb{\mathrm{comp}}(\mathcal{O}).
\end{equation}
By Lemma~\ref{lemma-bd-c},
the relation
\eqref{int-f-g-n} leads to
\[
\int_{\mathcal{O}}\rho(y)\int\sb{0}\sp{\omega}f(y,\omega-\tau)g(y,\tau)\tau^N\,d\tau\,dy=0,
\quad N\in\N,
\quad\forall\omega\in(0,\omega_0),
\quad\forall \rho\in C^\infty\sb{\mathrm{comp}}(\mathcal{O}).
\]
It follows that
\begin{equation}\label{f-g-zero}
f(y,\omega-\tau)g(y,\tau)=0,
\qquad \forall y\in \mathcal{O},
\qquad \forall\omega\in(0,\omega_0).
\end{equation}
Since we consider the case
$f\in \mathscr{D}(\R,L^2_{\mathrm{loc}}(\R^n))
\subset L^2_{\mathrm{loc}}(\R^n\times\R)$,
for a given open neighborhood $\mathcal{O}\ni x$
there is an open neighborhood
$\mathcal{O}_1\subset \mathcal{O}$,
an open interval
$\Omega_1\subset\big(\aa{f}(x)-\epsilon,\aa{f}(x)+\epsilon\big)$,
and $\delta>0$
such that
$\abs{f}\ge\delta$
almost everywhere on
$\mathcal{O}_1\times\Omega_1$.
(If not, then one would conclude that
$f=0$
almost everywhere
in
$\mathcal{O}\times\big(\aa{f}(x)-\epsilon,\aa{f}(x)+\epsilon\big)$,
contradicting the definition of
$\aa{f}(x)$.)
It follows from \eqref{f-g-zero}
that
$g(y,\omega-\tau)=0$
for all $y\in \mathcal{O}_1$,
$\omega\in(0,\omega_0)$,
$\tau\in\Omega_1$.
Therefore,
\[
g(y,\omega)
\equiv 0
\qquad
\mbox{almost everywhere in the rectangle}
\quad
(y,\omega)\in
\mathcal{O}_1\times(0,\omega_1),
\]
where $\omega_1:=\omega_0-\aa{f}(x)-\epsilon$.
Choosing $\epsilon=2^{-j}$, $j\in\N$, in the above construction,
we obtain a sequence
$\big(\omega_j\big)\sb{j\in\N}$
which converges to $\aa{f\astomega g}(x)-\aa{f}(x)$
and $\mathcal{O}_j\subset \mathbb{B}_{\epsilon_j}(x)$
such that
$
\big(\mathcal{O}_j\times(0,\omega_j)\big)
\cap\supp g=\emptyset.
$
(See Fig.~1).

\begin{figure}[htbp]
\begin{center}
\setlength{\unitlength}{1pt}
\begin{picture}(-145,130)(155,-70)
\put(80,35){$\supp g$}
\put(97,13.5){$g\equiv 0$}
\put(120,30){\circle*{4}}
\put(120,32){$\ \ \aa{f}^U(x)$}
\put(120,-30){\circle*{4}}
\put(120,-25){$\ \ \aa{f}(x)$}
\linethickness{0.2pt}
\put(120,-32){\line(0,1){4}}
\put(120,-37){\hskip -1pt $x$}
\put(35,5){\line(0,-1){35}}
\put(55,5){\line(0,-1){35}}
\put(35,5){\line(1,0){20}}
\multiput(10,4.5)(3,0){9}{.}
\put(40,-40){$\mathcal{O}_1$}
\put(-5,3){$\omega_1$}
\put(67,17){\line(0,-1){47}}
\put(82,17){\line(0,-1){47}}
\put(67,17){\line(1,0){15}}
\multiput(10,16.5)(3,0){19}{.}
\put(69,-40){$\mathcal{O}_2$}
\put(-5,15){$\omega_2$}
\put(92,25){\line(0,-1){55}}
\put(102,25){\line(0,-1){5}}
\put(102,-30){\line(0,1){40}}
\put(92,25){\line(1,0){10}}
\multiput(10,24.5)(3,0){27}{.}
\put(95,-40){$\mathcal{O}_3$}
\put(-5,24){$\omega_3$}
\put(10,-30){\vector(1,0){150}}
\put(10,-30){\vector(0,1){70}}
\linethickness{0.5pt}
\linethickness{1pt}
\put(120,30){\line(0,-1){60}}
\qbezier(120,-30)(145,-28)(150,-35)
\qbezier(10,-10)(40,30)(120,30)
\linethickness{0.1pt}
\multiput(36, 15)(1,-0){17}{\line(1,2){1}}
\multiput(42, 10)(1,-0){14}{\line(1,2){1}}
\multiput(45,  5)(1,-0){13}{\line(1,2){1}}
\multiput(46,  0)(1,-0){13}{\line(1,2){1}}
\multiput(45, -5)(1,-0){14}{\line(1,2){1}}
\multiput(42,-10)(1,-0){16}{\line(1,2){1}}
\multiput(38,-15)(1,-0){18}{\line(1,2){1}}
\multiput(40,-20)(1,-0){12}{\line(1,2){1}}
\put(50,-55){{\bf Fig. 1}\quad $\aa{f}$ and $\aa{f}^U$.}
\end{picture}
\end{center}
\end{figure}

It follows that
$
\aa{g}^U(x)\ge \aa{f\astomega g}(x)-\aa{f}(x),
$
and similarly
$
\bb{g}^L(x)\le \bb{f\astomega g}(x)-\bb{f}(x).
$
Since $x\in\R^n$ was arbitrary, this finishes the proof.
\end{proof}

\subsection{Proof of partial convolution theorem
for $f,\,g\in \mathscr{E}'(\R,L^2_{\mathrm{loc}}(\R^n))$}
\label{sect-f-g-dp}

\begin{lemma}\label{lemma-mollifiers}
Let $f,\,g\in \mathscr{E}'(\R,L^2_{\mathrm{loc}}(\R^n))$.
Let $\varphi\in\mathscr{D}(\R)$,
$\int\sb{\R}\varphi(\omega)\,d\omega=1$.
Then
$\supp f\astomega\varphi\to \supp f$
as $\supp\varphi\to\{0\}$,
and moreover, for each $x\in\R^n$,
\[
\aa{f\astomega\varphi}(x)\to \aa{f}(x),
\qquad
\aa{f\astomega\varphi}^U(x)\to \aa{f}^U(x)
\qquad
\mbox{as $\supp\varphi\to\{0\}$},
\]
\[
\bb{f\astomega\varphi}(x)\to \bb{f}(x),
\qquad
\bb{f\astomega\varphi}^L(x)\to \bb{f}^L(x)
\qquad
\mbox{as $\supp\varphi\to\{0\}$}.
\]
\end{lemma}


\begin{proof}
If $(x,\omega)\in\supp f$,
then there is an arbitrarily small open neighborhood
$\mathcal{O}\times\Omega$
of $(x,\omega)$
and the functions
$\psi\in L^2(\mathcal{O})$
and
$\theta\in \mathscr{D}(\Omega)$
such that
$\langle f,\psi\otimes\theta\rangle\ne 0$.
One has
$\theta\astomega\varphi
=\theta\ast\varphi
\mathop{\longrightarrow}\limits\sp{\mathscr{D}}
\theta$
as $\supp\varphi\to\{0\}$
(see \cite[Theorem 1.3.2]{MR717035});
then
\[
0\ne
\langle f,\psi\otimes\theta\rangle
=
\lim\sb{\supp\varphi\to \{0\}}
\langle f,\psi\otimes(\varphi\ast\theta)\rangle.
\]
Therefore,
one has
$
\langle f\astomega\varphi,\psi\otimes\theta\rangle
=\langle f,\psi\otimes(\varphi\ast\theta)\rangle
\ne 0$
for $\supp\varphi$ small enough.
For such $\varphi$,
one has
\[
\dist(\supp f,\supp f\astomega\varphi)
\le
\mathop{\rm diam}(\mathcal{O})
+
\mathop{\rm diam}(\Omega)
+
\mathop{\rm diam}(\supp \varphi)
.
\]
Since $\mathcal{O}$ and $\Omega$
are arbitrarily small, the conclusion follows.
\end{proof}

\begin{proof}
[Proof of Theorem~\ref{theorem-partial}]
We follow
the proof of \cite[Theorem 4.3.3]{MR717035}.
Let $0\le\varphi\in\mathscr{D}(\R)$ be such that
$\int_\R\varphi(\omega)\,d\omega=1$;
we apply
the version of Theorem~\ref{theorem-partial}
for $f,\,g\in
\mathscr{D}(\R,L^2\sb{\mathrm{loc}}(\R^n))$
(which we proved in Section~\ref{sect-f-g-d})
to $f\astomega \varphi$ and
$g\astomega\varphi$
to conclude that
\[
\bb{f\astomega\varphi}+\bb{g\astomega\varphi}^L
\le
\bb{(f\astomega\varphi)\astomega(g\astomega\varphi)}
=
\bb{(f\astomega g)\astomega(\varphi\astomega\varphi)}.
\]
Considering the limit $\supp\varphi\to\{0\}$
and applying Lemma~\ref{lemma-mollifiers},
we arrive at
\begin{eqnarray}\label{f-g-l-f-g}
\bb{f}(x)+\bb{g}^L(x)
\le
\bb{f\astomega g}(x),
\qquad
x\in\R^n.
\end{eqnarray}

\begin{lemma}\label{lemma-a-u-l-a}
Let $f\in 
\mathscr{E}'(\R,L^2_{\mathrm{loc}}(\R^n))$.
Then
$
(\aa{f}^U)^L=\aa{f}
$
and
$
(\bb{f}^L)^U=\bb{f}.
$
\end{lemma}

\begin{proof}
It is enough to prove the second statement.
Let us first prove it for $f$ measurable.
For $x\in\R^n\setminus\Sigma[f]$,
since $\Sigma[f]$ is closed
(see Lemma~\ref{lemma-sigma-closed}),
there is an open neighborhood
$\mathcal{O}\subset\R^n$,
$\mathcal{O}\ni x$,
such that $\mathcal{O}\cap\Sigma[f]=\emptyset$,
hence
\[
\bb{f}\at{\mathcal{O}}\equiv-\infty,
\qquad
\bb{f}^U\at{\mathcal{O}}\equiv-\infty,
\qquad
(\bb{f}^U)^L\at{\mathcal{O}}\equiv-\infty.
\]
Now let us consider $x\in\Sigma[f]\subset\R^n$.
For any $\epsilon>0$,
there is $\delta>0$,
$\mathcal{O}\subset\mathbb{B}_\epsilon(x)$,
and $\Omega\subset\big(\bb{f}(x)-\epsilon,\bb{f}(x)+\epsilon\big)$
such that $\abs{f}\ge\delta$
for almost all $(x,\omega)\in\mathcal{O}\times\Omega$
(otherwise, $f$ vanishes almost everywhere in an open neighborhood of
$(x,\bb{f}(x))$, hence $(x,\bb{f}(x))\not\in\supp f$,
which is in contradiction to the definition of $\bb{f}$).

For $f\in 
\mathscr{E}'(\R,L^2_{\mathrm{loc}}(\R^n))$,
we fix $\varphi\in C^\infty\sb{\mathrm{comp}}(\R)$,
$\varphi\ge 0$,
$0\in\supp\varphi$,
$\int_{\R}\varphi(\omega)\,d\omega=1$,
and consider $f\astomega\varphi$.
Since $f\astomega\varphi\in
\mathscr{D}(\R,L^2\sb{\mathrm{loc}}(\R^n))$
is measurable,
the first part of the proof applies,
showing that
\begin{eqnarray}\label{measurable-good}
(\bb{f\astomega\varphi}^L)^U(x)=\bb{f\astomega\varphi}(x),
\qquad
\forall x\in\R^n.
\end{eqnarray}
It remains to notice that
$
\bb{f}+\mathop{\mathrm{diam}}(\supp\varphi)
\ge
\bb{f\astomega\varphi}
\ge
\bb{f},
$
for all $x\in\R^n$,
with the last inequality due to
\eqref{f-g-l-f-g},
and to send $\supp\varphi\to\{0\}$;
then \eqref{measurable-good} turns into
$(\bb{f}^L)^U=\bb{f}$,
for all $x\in\R^n$.
\end{proof}

\begin{lemma}\label{lemma-reformulation}
Assume that
$f,\,g\in
\mathscr{E}'(\R,L^2\sb{\mathrm{loc}}(\R^n))$.
Then
\[
\big(\aa{f}+\aa{g}^U\big)^L
=
\aa{f\astomega g}
=
\big((\aa{f}+\aa{g})^U\big)^L,
\qquad
\big(\bb{f}+\bb{g}^L\big)^U
=
\bb{f\astomega g}
=
\big((\bb{f}+\bb{g})^L\big)^U.
\]
\end{lemma}

\begin{proof}
It is enough to prove the second statement.
From \eqref{f-g-l-f-g}, we conclude that
$
\big(\bb{f}+\bb{g}^L\big)^U\le \bb{f\astomega g},
$
while
\[
\bb{f\astomega g}\le \bb{f}+\bb{g}
\ \Rightarrow\ \bb{f\astomega g}^L\le \big(\bb{f}+\bb{g}\big)^L
\ \Rightarrow\ \big(\bb{f\astomega g}^L\big)^U
\le\big((\bb{f}+\bb{g})^L\big)^U.
\]
By Lemma~\ref{lemma-a-u-l-a},
$\big(\bb{f\astomega g}^L\big)^U=\bb{f\astomega g}$;
therefore, we conclude from the above relations that
\begin{eqnarray}\label{le-le}
\big(\bb{f}+\bb{g}^L\big)^U\le \bb{f\astomega g}
\le
\big((\bb{f}+\bb{g})^L\big)^U.
\end{eqnarray}
On the other hand, let us pick $x\in\R^n$;
there is a sequence
$x_j\to x$
such that
$\bb{g}(x_j)\to \bb{g}^L(x)$.
Then
$\lim\sup \bb{f}(x_j)\le \bb{f}(x)$,
hence we conclude that
\begin{eqnarray*}
\bb{f}(x)+\bb{g}^L(x)
=
\bb{f}(x)+\lim_{j\to\infty} \bb{g}(x_j)
\ge
\mathop{\lim\sup}_{j\to\infty}\big(\bb{f}(x_j)+\bb{g}(x_j)\big)
\ge
\mathop{\lim\inf}_{j\to\infty}
\big(\bb{f}(x_j)+\bb{g}(x_j)\big)
\ge
\big(\bb{f}+\bb{g}\big)^L(x),
\end{eqnarray*}
and then we conclude that
\begin{eqnarray}\label{ge-ge}
\big(\bb{f}+\bb{g}^L\big)^U(x)
\ge
\big((\bb{f}+\bb{g})^L\big)^U(x),
\qquad
\forall x\in\R^n.
\end{eqnarray}
Combining \eqref{le-le} and \eqref{ge-ge},
we arrive at
$
\ \big(\bb{f}(x)+\bb{g}^L\big)^U(x)
=
\bb{f\astomega g}(x)
=
\big(\big(\bb{f}+\bb{g}\big)^L\big)^U(x),
$
\ for all $x\in\R^n$.
\end{proof}

This completes the proof of
Theorem~\ref{theorem-partial}.
\end{proof}



\section{Compact spectrum solutions to the nonlinear Klein--Gordon equation}
\label{sect-kg}

Let us first recall the unique continuation property ({\bf UCP})
for the Laplace operator.

\begin{theorem}
[Unique continuation property  for the Laplace operator]
\label{theorem-ucp}
Let $n\ge 1$.
Assume that
$u\in H^1(\R^n)$ satisfies the relation
\begin{align}\label{ucp-laplace}
\abs{\Delta u}\le\abs{V u}
\end{align}
almost everywhere in a connected open domain
$\mathcal{O}\subset\R^n$,
with $V\in L^p\sb{\mathrm{loc}}(\R^n)$,
with $p\ge n/2$, $n\ge 2$,
and with $p=1$ for $n=1$.
If $u$ vanishes almost everywhere in an open subset
$\mathcal{O}_0\subset \mathcal{O}$,
then it vanishes almost everywhere in $\mathcal{O}$.
\end{theorem}

Wolff \cite[Theorem 3]{MR1159832}
proved the unique continuation in $\R^n$, $n\ge 3$,
with
$V\in L^{p}\sb{\mathrm{loc}}(\R^n)$,
where $p\ge n/2$ if $n\ge 5$, $p>2$ if $n=4$, and $p\ge 2$ if $n=3$.
The optimal unique continuation results for
\eqref{ucp-laplace}
were obtained in \cite{MR1809741};
in particular, it follows that the strong unique continuation property
holds for
$V\in L^{n/2}\sb{\mathrm{loc}}(\R^n)$
for any $n\ge 2$
(the sufficient conditions on $V$
in \cite[Theorem 1.1]{MR1809741} are slightly weaker).
We also mention that the
sufficient condition $V\in L^1\sb{\mathrm{loc}}(\R)$
in the one-dimensional case
is a consequence of the two-dimensional result
when considering functions with a trivial dependence on $x_2$.
Or, arguing directly,
one could assume that $I\subset\R$ is an open interval
and $u\in H^1(\R)$ vanishes in an open neighborhood
$\mathcal{O}_0\subset I$
and satisfies
$\abs{u''}\le \abs{V u}$ almost everywhere in $I$,
with some $V\in L^1\sb{\mathrm{loc}}(\R)$.
Let us show that $u\at{I}\equiv 0$.
Without loss of generality, we may assume that $I$ is bounded.
Given $x_0\in\mathcal{O}_0$,
then, for any $x\in I$,
one has
$\abs{u'(x)}\le\int_{x_0}^x\abs{V(y)u(y)}\,dy
\le\sup\sb{y\in[x_0,x]}\abs{u(y)}\norm{V}\sb{L^1(I)}$,
hence
\[
\abs{u(x)}\le\int_{x_0}^x\abs{u'(y)}\,dy
\le\sup\sb{y\in[x_0,x]}\abs{u(y)}
\abs{x-x_0}\norm{V}\sb{L^1(I)},
\]
showing that $u(x)=0$ as long as
$x\in I$ is close enough to $x_0$ so that
$\abs{x-x_0}<1/\norm{V}\sb{L^1(I)}$.
It follows that $u\equiv 0$ in $I$.

\begin{remark}\label{remark-ucp}
For $n\le 2$,
the Sobolev embedding
gives $u\in H^1(\R^n)\subset L^q(\R^n)$
for any $2\le q<\infty$
(including $q=\infty$ when $n=1$),
hence
$V(x):=\alpha(\abs{u(x)}^2)$
with $\alpha(\tau)$ from \eqref{kappa-such}
satisfies
$V\in L^p\sb{\mathrm{loc}}(\R^n)$
for any $1<p<\infty$.
Therefore, for $n\le 2$,
the unique continuation takes place for any $\kappa>0$.

For $n\ge 3$,
by the Sobolev embeddings,
$
u\in
H^1(\R^n)\subset L^{2^*}(\R^n),
$
with
$2^*=2n/(n-2)$.
Then $V(x):=\alpha(\abs{u(x)}^2)$
satisfies
\[
V\in L^p\sb{\mathrm{loc}}(\R^n),
\qquad
\mbox{with}
\quad
p=\frac{2^*}{2\kappa}
=\frac{n}{(n-2)\kappa}.
\]
For the unique continuation to take place,
we need the relation
$
p=\frac{n}{(n-2)\kappa}\ge \frac{n}{2}
$,
so for $n\ge 3$ we need
$\kappa\le 2/(n-2)$.
\end{remark}



Now we recall the local well-posedness results for the Klein--Gordon equation.


\begin{theorem}[NLKG global well-posedness
\mbox{\cite[Proposition 2.1]{MR843591}}]
\label{theorem-kg-gwp}
Let $n\in\N$, $m>0$.
Let $f\in C^1(\C,\C)$ with $f(0)=0$
and $f(e^{\jj s}u)=e^{\jj s}f(u)$,
$\forall u\in\C$, $\forall s\in\R$;
$F(u)=\int_0^{\abs{u}} f(v)\,dv$,
$u\in\C$.
Assume that there are $c_0>0$ and $c_1>0$
such that
\[
F(u)\ge -c_0\abs{u},
\quad
\abs{f'(u)}
\le c_1(1+\abs{u}^{p-1}),
\quad \forall u\in\C,
\]
with some $p\in(1,+\infty)$ if $n\le 2$;
$p\in(1,1+4/(n-2))$ if $n\ge 3$.
Then there is a unique, strongly continuous solution
$u\in C(\R,H^1(\R^n,\C))$,
$\forall t\in\R$,
to the Cauchy problem
\begin{eqnarray}\label{nlkg-psi}
-\p_t^2 u=-\Delta u+m^2 u+f(u),
\quad
u(x,t)\in\C,
\quad
x\in\R^n;
\qquad
(u,\dot u)\at{t=0}\in H^1(\R^n,\C)\times L^2(\R^n,\C).
\end{eqnarray}
Its energy is conserved:
$E(u(t))=E(u(0))$
\ for all \ $t\in\R$.
\end{theorem}

\begin{assumption}\label{ass-kappa}
$f(u)=\alpha(\abs{u}^2)u$,
with
$\alpha\in C\big(\overline{\R\sb{+}}\big)$,
$\alpha(0)=0$,
and there is $C<\infty$
such that
\begin{eqnarray}\label{kappa-such}
\abs{\alpha(\tau)}\le C\langle\tau\rangle^\kappa,
\qquad
\forall\tau\ge 0,
\qquad
\mbox{with $\kappa$ satisfying}
\quad
\begin{cases}
\kappa>0,& n\le 2;
\\
0<\kappa\le 2/(n-2),&
n\ge 3.
\end{cases}
\end{eqnarray}
\end{assumption}
We note that the restriction on $\kappa$
is such that the unique continuation property
from Theorem~\ref{theorem-ucp}
applies to $V(x)=\alpha(\abs{u}^2)$
with $u\in H^1(\R^n)$
(see Remark~\ref{remark-ucp}).
We also note that the well-posedness result from Theorem~\ref{theorem-kg-gwp}
applies
if e.g.
$\alpha(\tau)=C\tau^\kappa$, $\tau\ge 0$
(or if $\alpha(\tau)$ is a polynomial of degree $\kappa\in\N$),
with $\kappa>0$ for $n\le 2$, $0<\kappa<2/(n-2)$
if $n\ge 3$
(cf. \eqref{kappa-such}).

\medskip

We will be able to consider not only polynomial nonlinearities,
but also certain algebraic nonlinearities.

\begin{assumption}\label{ass-algebraic}
Assume that $\alpha\in C\big(\overline{\R\sb{+}}\big)$ is a non-constant algebraic function,
so that there is $J\in\N$
and polynomials
$M_j(\tau)$, $0\le j\le J$,
with $M_J(\tau)\not\equiv 0$,
such that
$w(\tau):=\tau \alpha(\tau)$
satisfies the relation
$\mathcal{M}(\tau,w(\tau))=0$, $\forall\tau\ge 0$,
where
\begin{eqnarray}\label{def-poly-Q}
\mathcal{M}(\tau,w(\tau)):=
\tau^J\sum_{j=0}^J P_j(\tau)\alpha(\tau)^j
=\sum_{j=0}^J M_j(\tau)(\tau \alpha(\tau))^j
=\sum_{j=0}^J M_j(\tau)w(\tau)^j,
\qquad
\forall\tau\ge 0.
\end{eqnarray}
Moreover, assume that
\[
\deg M_0>\deg M_j+j,\qquad
\forall j,
\quad 1\le j\le J.
\]
If $n\ge 3$,
additionally assume that
\begin{eqnarray}\label{q-k}
\deg M_j+(\kappa+1)j\le n/(n-2),
\qquad
\forall j,\quad 0\le j\le J.
\end{eqnarray}
\end{assumption}


\begin{example}
Assume that
$\alpha$ is a polynomial:
$\alpha(\tau)=\sum_{j=0}^\kappa \alpha_j\tau^j$,
with $\kappa\in\N$,
$\kappa\le 2/(n-2)$ if $3\le n\le 4$,
and $\alpha_\kappa\ne 0$.
Let
$M_0(\tau)=-\tau \alpha(\tau)$
and
$M_1(\tau)=1$,
so that
$\deg M_0=\kappa+1$
and $\deg M_1=0$.
Then
\[
\mathcal{M}(\tau,\tau \alpha(\tau))
=
M_0(\tau)
+
M_1(\tau)\tau \alpha(\tau)
=
-\tau \alpha(\tau)\cdot 1+1\cdot\tau \alpha(\tau)=0,
\qquad
\forall\tau\in\R.
\]
One can see that Assumption~\ref{ass-algebraic} is satisfied
(including the requirement \eqref{q-k} when $n\ge 3$).
\end{example}

\begin{example}
Assume that $\alpha(\tau)=A(\tau)^{1/N}$,
with $N\in\N$, $N\ge 2$,
and with
$A(\tau)=\sum\sb{j=0}^a A_j\tau^j$
a polynomial of degree $a=\deg A\ge 1$;
if $N$ is even, we additionally assume that
$A(\tau)\ge 0$ for $\tau\ge 0$.
Let
$M_0(\tau)=-\tau^N A(\tau)$,
and $M_N(\tau)=1$;
$\deg M_0=a+N$ and $\deg M_N=0$.
Then
\[
\mathcal{M}(\tau,\tau \alpha(\tau))
=M_0(\tau)+M_N(\tau)(\tau \alpha(\tau))^N
=
-\tau^N A(\tau)\cdot 1+1\cdot(\tau \alpha(\tau))^N=0,
\qquad
\forall\tau\ge 0.
\]
If $n\le 2$,
Assumption~\ref{ass-algebraic} is satisfied.
If $n\ge 3$,
we additionally need
$0<\kappa=a/N\le 2/(n-2)$
and \eqref{q-k} to be satisfied;
in this case,
there are nontrivial examples only when $a=1$, $N=2$, $n=3$.
\end{example}

\begin{example}
Consider
$\alpha(\tau)=A(\tau)/B(\tau)$,
with $A,\,B$ polynomials of degrees $a=\deg A\ge 0$ and $b=\deg B\ge 1$;
$B(\tau)\ne 0$ for $\tau\ge 0$.
Let
$M_0(\tau)=-\tau A(\tau)$,
$M_1(\tau)=B(\tau)$.
Then
\[
\mathcal{M}(\tau,\tau \alpha(\tau))
=-\tau A(\tau)\cdot 1+B(\tau)\cdot\tau \alpha(\tau)=0,
\qquad
\forall\tau\ge 0.
\]
For Assumption~\ref{ass-algebraic} to be satisfied,
we need
$
\deg M_0=a+1>\deg M_1+1=b+1,
$
so for $n\le 2$ one only needs $a>b$.
If $n\ge 3$, one additionally needs
$0<\kappa=a-b\le 2/(n-2)$
and \eqref{q-k};
since $a>b\ge 1$,
there are nontrivial examples only when $n=3$, $a=2$, $b=1$.
\end{example}

Now we can formulate and prove our main result:
under rather generic assumptions
the only type of solutions with compact time spectrum
is the one-frequency solitary waves.

\begin{theorem}
\label{theorem-one-frequency}
Let $n\in\N$, $m>0$.
Let $f(\upsi)=\alpha(\abs{\upsi}^2)\upsi$
be such that
both
Assumption~\ref{ass-kappa}
and 
Assumption~\ref{ass-algebraic}
are satisfied.
Assume that $\upsi\in L^\infty(\R,H^1(\R^n))$
is a solution to the nonlinear Schr\"odinger or Klein--Gordon equation
\eqref{eqns}.
If there is a finite interval $I\subset\R$
such that $\supp\tilde\upsi\subset\R^n\times I$,
with $\tilde\upsi(x,\omega)$ the Fourier transform of $\upsi$
with respect to time,
then
\[
\upsi(x,t)=\phi_0(x)e^{-\jj\omega_0 t},
\qquad
\mbox{
with some \ $\phi_0\in H^1(\R^n,\C)$ \ and \ $\omega_0\in\R$}.
\]
\end{theorem}

Note that, in particular, the above theorem applies to finite energy solutions
to the nonlinear Klein--Gordon equation
from Theorem~\ref{theorem-kg-gwp}.

\begin{proof}
The proof for the nonlinear Schr\"odinger equation
and the nonlinear Klein--Gordon equation is the same;
for definiteness, we consider the latter case.
Assume that $\upsi\in L^\infty(\R,H^1(\R^n,\C))$
is a solution to
\eqref{nlkg-psi}
with compact time spectrum,
so that the Fourier transform of $\upsi$ in time,
\[
\tilde\upsi(x,\omega)
=
\int\sb{\R}\upsi(x,t)e^{\jj\omega t}\,dt,
\qquad
\tilde\upsi\in\mathscr{E}'(\R,H^1(\R^n,\C)),
\]
satisfies
$\supp\tilde\upsi\subset\R^n\times [a,b]$,
with some $a,\,b\in\R$, $a<b$.
We denote
\begin{eqnarray}\label{def-sigma-psi}
\Sigma:=\Sigma[\tilde\upsi]
=\big\{
x\in\R^n\sothat
(\{x\}\times\R)\cap\supp\upsi\ne\emptyset
\big\}
=\big\{
x\in\R^n\sothat
\big(\{x\}\times\R\big)\cap\supp\tilde\upsi\ne\emptyset
\big\}
\end{eqnarray}
to be the projection of the support of $\upsi$
onto $\R^n$.
Then,
since
$\supp\tilde\upsi\subset\R^n\times[a,b]$,
\[
\bb{\tilde\upsi}\at{\Sigma}\ge a
\quad\Rightarrow\quad
\bb{\tilde\upsi}^L\at{\Sigma\setminus\p\Sigma}\ge a;
\qquad
\aa{\tilde\upsi}\at{\Sigma}\le b
\quad\Rightarrow\quad
\aa{\tilde\upsi}^U\at{\Sigma\setminus\p\Sigma}\le b.
\]

\begin{lemma}\label{lemma-b-l-a}
$\alpha(\abs{\upsi(x,t)}^2)$ and $\abs{\upsi(x,t)}$
do not depend on time,
and moreover
\[
\bb{\tilde\upsi}^L=\aa{\tilde\upsi},
\qquad
\bb{\tilde\upsi}=\aa{\tilde\upsi}^U,
\qquad
\forall x\in\Sigma.
\]
\end{lemma}

\begin{proof}
The Sobolev embedding
leads to
\begin{eqnarray}\label{psi-q}
\upsi
\in L^\infty(\R,H^1(\R^n,\C))
\subset L^\infty(\R,L^q\sb{\mathrm{loc}}(\R^n,\C)),
\qquad
\begin{cases}
1\le q<\infty,& n\le 2;
\\
1\le q\le 2n/(n-2),&n\ge 3.
\end{cases}
\end{eqnarray}
The inclusion \eqref{psi-q}
together with \eqref{kappa-such}
lead to
\begin{eqnarray}\label{g-is-lq}
v\in L^\infty\big(\R,L^{q/(2\kappa)}\sb{\mathrm{loc}}(\R^n,\R)\big),
\qquad
v(x,t):=\alpha(\abs{\upsi(x,t)}^2).
\end{eqnarray}
By \eqref{nlkg-psi} and \eqref{g-is-lq},
\begin{eqnarray}\label{w-is-l-q}
(\p_t^2-\Delta+m^2)\upsi
=-\alpha(\abs{\upsi}^2)\upsi\in L^\infty\big(\R,L^{q/(2\kappa+1)}\sb{\mathrm{loc}}(\R^n,\C)
\big).
\end{eqnarray}


Applying the Fourier transform to \eqref{w-is-l-q}
and denoting by
$\tilde v(x,\omega)$
the Fourier transform of $v(x,t):=\alpha(\abs{\upsi(x,t)}^2)$
in time, one has
\begin{eqnarray}\label{bar-psi-psi}
(m^2-\omega^2-\Delta)\tilde\upsi
=-\tilde v\astomega\tilde\upsi.
\end{eqnarray}

\noindent
$\bullet$
Let us
consider the case when $\alpha(\tau)$
is a polynomial of degree $\kappa=\deg \alpha\ge 1$,
with either $n\le 2$, $\kappa\in\N$;
or $n=3$, $\kappa=1,\,2$;
or $n=4$, $\kappa=1$.
Applying Theorem~\ref{theorem-partial}
to the right-hand side of the above relation,
we arrive at
\[
\bb{(m^2-\omega^2-\Delta)\tilde\upsi}(x)
=
\bb{\tilde v\astomega\tilde\upsi}(x)
\ge
\bb{\tilde v}^L(x)+\bb{\tilde\upsi}(x),
\qquad
\forall x\in\Sigma.
\]
Due to the inclusion
$\supp \Delta \tilde\upsi\subset\supp\tilde\upsi$,
the above yields
$
\bb{\tilde\upsi}
\ge
\bb{\tilde v}^L+\bb{\tilde\upsi}
$
for all $x\in\Sigma$,
hence
$\bb{\tilde v}^L\le 0$
and therefore
\begin{eqnarray}\label{b-l-negative}
\bb{\tilde v}(x)\le 0,
\qquad
\forall x\in\Sigma.
\end{eqnarray}
Similarly,
$
\aa{\tilde v}\ge 0,
$
for all $x\in\Sigma$;
thus,
$\supp\tilde v\subset\R^n\times\{0\}$,
and we conclude that
\begin{eqnarray}\label{g-powers}
v(x,t)=\sum_{j\in\N\sb{0}}t^j v_j(x).
\end{eqnarray}
Above,
in agreement with the general theory of distributions
\cite{MR717035},
the summation in $j\in\N\sb{0}$ is locally finite
(there are finitely many terms for $x\in K$
for each compact subset $K\subset\R^n$)
(cf. \cite[Theorem 2.3.5]{MR717035}).
The terms with derivatives of $\delta(\omega)$
do not appear since this would lead to $v(x,t)$ growing in time,
contradicting \eqref{g-is-lq}.
This implies that in \eqref{g-powers} the only nonzero term is the one
with $j=0$.
Thus, $V(x):=v(x,t)=\alpha(\abs{\upsi(x,t)}^2)$ does not depend on time.
Since $\alpha(\tau)$ is a nonconstant algebraic function,
$\abs{\upsi(x,t)}^2$ also does not depend on time:
\begin{eqnarray}\label{zero-is-b}
\supp\widetilde{\abs{\upsi}^2}\subset\R^n\times\{0\}.
\end{eqnarray}
Using the above relation and
applying Theorem~\ref{theorem-partial}
to
$\widetilde{\abs{\upsi}^2}=
\tilde\upsi\sp\sharp\astomega\tilde\upsi$,
where
$\tilde\upsi\sp\sharp=\widetilde{\bar\upsi}$
(see \eqref{def-sharp}),
we conclude that
\[
0
=\bb{\widetilde{\abs{\upsi}^2}}(x)
\ge
\bb{\tilde\upsi}^L(x)+\bb{\tilde\upsi\sp\sharp}(x)
=\bb{\tilde\upsi}^L(x)-\aa{\tilde\upsi}(x),
\qquad
\forall x\in\Sigma.
\]
Thus,
$\bb{\tilde\upsi}^L\le\aa{\tilde\upsi}$
for all $x\in\Sigma$.
On the other hand, by Lemma~\ref{lemma-a-a-b},
$\bb{\tilde\upsi}^L\ge \aa{\tilde\upsi}$
for all $x\in\Sigma$.
We conclude that
\begin{eqnarray}\label{a-is-b}
\bb{\tilde\upsi}^L=\aa{\tilde\upsi}
\quad
\mbox{and similarly}
\quad
\bb{\tilde\upsi}=\aa{\tilde\upsi}^U,
\qquad
\forall x\in\Sigma.
\end{eqnarray}

\noindent
$\bullet$
Let us consider the case when
$\alpha(\tau)$ is an algebraic function satisfying Assumption~\ref{ass-algebraic}.
Multiplying \eqref{w-is-l-q}
by $\bar\upsi$,
we have:
\begin{eqnarray}\label{psi-psi}
\bar\upsi(m^2+\p_t^2-\Delta)\upsi
=-\abs{\upsi}^2 \alpha(\abs{\upsi}^2)
\in L^\infty\big(\R,L\sb{\mathrm{loc}}^{q/(2\kappa+2)}(\R^n,\C)\big),
\end{eqnarray}
with $1\le q<\infty$ if $n\le 2$
and $1\le q\le 2n/(n-2)$ if $n\ge 3$.
Let $\mathcal{M}$ be as in \eqref{def-poly-Q}.
Applying $\mathcal{M}(\abs{\upsi}^2,\cdot)$
to both sides of the relation \eqref{psi-psi}
leads to
\begin{eqnarray}\label{leads-to}
0=\mathcal{M}
\big(
\abs{\upsi}^2,
\abs{\upsi}^2 \alpha(\abs{\upsi}^2)
\big)
=
\mathcal{M}
\big(
\abs{\upsi}^2,
-\bar\upsi(m^2-\omega^2-\Delta)\upsi
\big)
=
\sum\sb{j=0}\sp J
M_j(\abs{\upsi}^2)
\big(
-\bar\upsi(m^2+\p_t^2-\Delta)\upsi
\big)^j.
\end{eqnarray}
We need to make sure that the right-hand side is a well-defined distribution.
Taking into account
\eqref{psi-q} and \eqref{w-is-l-q},
we conclude that
all the terms in the right-hand side are in
$\mathscr{E}'\big(\R,L^1\sb{\mathrm{loc}}(\R^n)\big)$
as long as
in \eqref{psi-q}
one can take $q\ge 1$
such that
\[
\frac{2\deg M_j}{q}
+\frac{2\kappa+2}{q}j
\le 1,
\qquad
\forall j,\quad 0\le j\le J.
\]
For $n\le 2$, we can satisfy the above
by taking $1\le q<\infty$ arbitrarily large;
for $n\ge 3$, the above is satisfied with $q=2n/(n-2)$
due to
the inequality \eqref{q-k} in Assumption~\ref{ass-algebraic}.

We note that
$\tilde\upsi\sp\sharp\astomega\tilde\upsi=\widetilde{\abs{\upsi}^2}$
and that
$
\supp
(\tilde\upsi\sp\sharp\astomega(m^2-\omega^2-\Delta)\tilde\upsi)
\subset
\supp(\tilde\upsi\sp\sharp\astomega\tilde\upsi),
$
hence
\begin{eqnarray}\label{small-supp}
\bb{\tilde\upsi\sp\sharp\astomega(m^2-\omega^2-\Delta)\tilde\upsi}
\le
\bb{\tilde\upsi\sp\sharp\astomega\tilde\upsi},
\qquad
\forall x\in\Sigma.
\end{eqnarray}
Now we apply Theorem~\ref{theorem-partial} to the Fourier transform
(in time)
of the relation \eqref{leads-to}
and use Assumption~\ref{ass-algebraic},
arriving at
\[
\bb{\tilde\upsi\sp\sharp\astomega(m^2-\omega^2-\Delta)\tilde\upsi}^L
\le 0;
\]
then
$\bb{\tilde\upsi\sp\sharp\astomega(m^2-\omega^2-\Delta)\tilde\upsi}\le 0$,
and similarly
$\aa{\tilde\upsi\sp\sharp\astomega(m^2-\omega^2-\Delta)\tilde\upsi}\ge 0$.
It follows that
\[
\supp\tilde\upsi\sp\sharp\astomega(m^2-\omega^2-\Delta)\tilde\upsi\subset\R^n\times\{0\},
\]
hence, by the argument after \eqref{g-powers},
$\abs{\upsi}^2 \alpha(\abs{\upsi}^2)$ is time-independent,
and so is $\abs{\upsi}^2$
(we note that
$\tau \alpha(\tau)$ is a nonconstant function of $\tau$:
indeed, if we had
$0=\mathcal{M}(\tau,\tau \alpha(\tau))=-C+\tau \alpha(\tau)$,
then
$M_0(\tau)=C$ and $M_1(\tau)=1$,
not satisfying Assumption~\ref{ass-algebraic}).
Therefore,
we again arrive at
\eqref{zero-is-b} and then \eqref{a-is-b} follows.
\end{proof}

By Lemma~\ref{lemma-b-l-a},
\begin{eqnarray}\label{g-tilde-good}
V(x):=v(x,t)
=
\alpha(\abs{\upsi(x,t)}^2)
\quad
\mbox{does not depend on time;}
\qquad
\tilde v(x,\omega)=2\pi\delta(\omega)V(x).
\end{eqnarray}
Due to \eqref{g-tilde-good},
equation \eqref{nlkg-psi}
takes the form
\begin{eqnarray}\label{psi-tilde-ucp}
\Delta\tilde\upsi
=
m^2\tilde\upsi
-\omega^2\tilde\upsi
+V(x)\tilde\upsi.
\end{eqnarray}
By \eqref{psi-q},
$\abs{\upsi}^2
\in 
L^\infty\big(\R,L^{q/2}\sb{\mathrm{loc}}(\R^n,\R)\big)$,
$\widetilde{\abs{\upsi}^2}
\in
\mathscr{E}'\big(\R,L^{q/2}\sb{\mathrm{loc}}(\R^n,\R)\big)
$
(we took into account the assumption that the spectrum of $\upsi$ is compact),
with any $q\ge 1$ if $n\le 2$
and $1\le q\le 2n/(n-2)$ if $n\ge 3$.
Then, according to the assumption \eqref{kappa-such},
\begin{eqnarray}\label{v-good}
V(x)=\alpha(\abs{\upsi(x,t)}^2)
\quad
\mbox{satisfies}
\quad
V\in L^{q/(2\kappa)}\sb{\mathrm{loc}}(\R^n,\R),
\end{eqnarray}
with any $q\ge 2\kappa$ for $n\le 2$
and $q=2n/(n-2)$ for $n\ge 3$.
Due to the requirement \eqref{kappa-such} on $\kappa$,
the function $V(x)$ satisfies conditions needed for the unique continuation property
(see Theorem~\ref{theorem-ucp}
and Remark~\ref{remark-ucp}).

Let us show that $\Sigma[\tilde\upsi]$
defined in \eqref{def-sigma-psi}
has to be the whole space.

\begin{lemma}\label{lemma-sigma-is-all}
If $\upsi$ is not identically zero, then
$\Sigma[\tilde\upsi]=\R^n$.
\end{lemma}

\begin{proof}
Assume that, on the contrary, $\Sigma[\tilde\upsi]\subsetneq\R^n$;
since $\Sigma[\tilde\upsi]$ is closed,
there is a nonempty connected open subset
$\mathcal{O}\subset\R^n$
such that $\mathcal{O}\cap\Sigma[\tilde\upsi]=\emptyset$.
Let $\Omega\subset\R$ be an open interval;
since $\mathcal{O}\cap\Sigma[\tilde\upsi]=\emptyset$,
one has
$\big(\mathcal{O}\times\Omega\big)\cap\supp\tilde\upsi=\emptyset$.
Since $V$ satisfies the assumptions of  Theorem~\ref{theorem-ucp},
we apply the unique continuation property
to an $L^2$-function
$\tilde\upsi$
(valued in $\mathscr{D}'(\Omega)$)
which solves
\eqref{psi-tilde-ucp},
concluding that
\[
\big(\R^n\times\Omega\big)\cap\supp\tilde\upsi=\emptyset.
\]
When applying the unique continuation property to
\eqref{psi-tilde-ucp},
we need to mention that
the multiplication by $\omega^2$
is a continuous automorphism in
$\mathscr{D}(\Omega)=C^\infty_{\mathrm{comp}}(\Omega)$
(in the Fr\'echet topology based on sup-norms in $C^k_{\mathrm{comp}}(\Omega)$, $k\ge 0$),
and hence also in $\mathscr{D}'(\Omega)$.
\end{proof}

\begin{lemma}\label{lemma-one-omega}
There is $\omega_0\in\R$ such that
$\supp\tilde\upsi\subset\R^n\times\{\omega_0\}$.
\end{lemma}

\begin{proof}
Pick $x_1\in\Sigma[\tilde\upsi]=\R^n$.
Denote $\omega_1=\bb{\tilde\upsi}(x_1)$.
We will show that for any open neighborhood $\Omega\subset\R$,
$\inf\Omega>\omega_1$,
one has
$\big(\R^n\times\Omega\big)\cap\supp\tilde\upsi=\emptyset$.

Since $\bb{\tilde\upsi}$ is upper semicontinuous,
for any $\epsilon>0$,
which we choose to be $\epsilon:=\dist(\omega_1,\Omega)/2>0$,
there is an open neighborhood $\mathcal{O}\subset\R^n$,
$\mathcal{O}\ni x_1$,
such that
$\bb{\tilde\upsi}\at{\mathcal{O}}<\omega_1+\epsilon$.
Let
$\varphi\in C^\infty_{\mathrm{comp}}(\R,\R)$,
$\supp\varphi\subset\Omega$.
Using the unique continuation property
exactly as in Lemma~\ref{lemma-sigma-is-all}, we conclude that
$\big(\mathcal{O}\times\Omega\big)\cap\supp\tilde\upsi=\emptyset$
implies that
$\big(\R^n\times\Omega\big)\cap\supp\tilde\upsi=\emptyset$.
Since the choice of $x_1\in\R^n$ was arbitrary, we conclude that
\[
\supp\tilde\upsi\subset\R^n\times(-\infty,\inf \bb{\tilde\upsi}].
\]
Similarly one proves that
\[
\supp\tilde\upsi\subset\R^n\times[\sup\aa{\tilde\upsi},+\infty).
\]
By Lemma~\ref{lemma-b-l-a},
$\bb{\tilde\upsi}^L=\aa{\tilde\upsi}$;
it follows that $\inf \bb{\tilde\upsi}=\sup\aa{\tilde\upsi}=:\omega_0$,
and therefore
$\supp\tilde\upsi\subset\R^n\times\{\omega_0\}$.
\end{proof}

By Lemma~\ref{lemma-one-omega},
\begin{eqnarray}\label{phi-j}
\upsi(x,t)= e^{-\jj\omega_0 t}\sum\sb{j\in\N\sb{0}}\phi_j(x)t^j.
\end{eqnarray}
By the above arguments,
the summation in \eqref{phi-j} is locally finite for $x\in K$,
for each compact subset $K\subset\R^n$.
Since $\upsi\in L^\infty(\R,H^1(\R^n,\C))$,
we conclude that
in \eqref{phi-j}
the terms with $j\ge 1$ are absent;
thus,
$\upsi(x,t)=\phi_0(x)e^{-\jj\omega_0 t}$,
with $\phi_0\in H^1(\R^n,\C)$.
This concludes the proof of Theorem~\ref{theorem-one-frequency}.
\end{proof}



\bibliographystyle{sima-doi}
\bibliography{bibcomech}
\end{document}